\documentclass[12pt,a4paper]{article}
\usepackage[utf8]{inputenc}
\usepackage[bitstream-charter]{mathdesign}
\usepackage{amsmath,graphicx}
\usepackage{bbm,float,longtable,authblk}
\usepackage{apacite,color,soul,xcolor,ntheorem,lipsum}
\usepackage[round]{natbib}
\usepackage[left=2.5cm,right=2.5cm,top=2.5cm,bottom=2.5cm]{geometry}

\usepackage{bm}
\usepackage{hyperref}
\theorembodyfont{\bfseries\itshape}

\newtheorem{rem}{Remark} 
\newtheorem{lemma}{Lemma}
\newtheorem{result}{Result}

\title{Some Permutationllay Symmetric Multiple Hypotheses Testing Rules Under Dependent Set up}

\author{Anupam Kundu\thanks{Department of Statistics, Texas A\& M University; Email:\url{akundu@stat.tamu.edu}} \and Subir Kumar Bhandari\thanks{Interdisciplinary Statistical Research Unit, Indian Statistical Institute, Kolkata; Email: \url{subir@isical.ac.in}}}

\date{}

\begin{document}
\maketitle

\begin{abstract}
In this paper, our interest is in the problem of simultaneous hypothesis testing when the test statistics corresponding to the individual hypotheses are possibly correlated. Specifically, we consider the case when the test statistics together have a multivariate normal distribution (with equal correlation between each pair) with an unknown mean vector and our goal is to decide which components of the mean vector are zero and which are non-zero. This problem was taken up earlier in \cite{bogdan2011asymptotic} for the case when the test statistics are independent normals. Asymptotic optimality in a Bayesian decision theoretic sense was studied in this context, the optimal precodures were characterized and optimality of some well-known procedures were thereby established. The case under dependence was left as a challenging open problem. We have studied the problem both theoretically and through extensive simulations and have given some permutation invariant rules. Though in \cite{bogdan2011asymptotic}, the asymptotic derivations were done in the context of sparsity of the non-zero means, our result does not require the assumption of sparsity and holds under a more general setup.
\end{abstract}



\section{Introduction}
Multiple hypothesis testing has emerged as a very important topic of research in the last twenty years. The biggest impetus to such work came from the necessity to analyze and draw inference on data sets involving a large number of parameters. Such data sets occur, e.g., in the fields of biology, astronomy, economics, just to name a few. Needless to say that the goal of simultaneous testing, or for that matter, simultaneous inference in general, is to ensure a good performance of the overall inference, rather than ensuring good performance for the individual inference problems. 
\par

Over the years, various performance evaluation criteria have been developed to quantify the overall error in a simultaneous testing procedure. The most classical measure of this kind is the Family Wise Error Rate (FWER).Well-known procedures that control the FWER are the Bonferroni procedure and its improvements, see for example, \citet{holm1979simple}, \citet{simes1986improved}, \citet{hommel1988stagewise} and \citet{benjamini1995controlling}. A nice historical account of the early works in this area can be found in \citet{hochberg1987multiple}. A great leap forward in the field of simultaneous inference was made through the introduction of the concept of false discovery rate (FDR) and a procedure called the Benjamini-Hochberg procedure. These appeared in the seminal work by Benjamini and Hochberg in 1995. FDR is obviously the more appropriate error to control in a large scale simultaneous testing problems compared to the FWER, since trying to control the probability of a single erroneous rejection seems too stringent a requirement in such cases. See, for example, \cite{benjamini1999step}, \cite{sarkar2007}, \cite{storey2002direct} and \cite{storey2004strong} for further details. An excellent account on the literature on FDR can be found in \cite{Sarkar2008} and \cite{efron2012large}\par

The degree of "surprise" required in the observed data to declare a particular hypothesis to be false in a multiple hypothesis testing context should be more than what would be required to reject a hypothesis in an individual testing problem. The examples given above of multiple testing procedures belong to the frequentist domain.  
For the Bayesian, it is intuitive to reject a hypothesis when it is less likely a posteriori.  The articles \cite{SCOTT20062144} and \cite{scott2010bayes} beautifully explain this insight and explicitly demonstrate multiplicity adjustment in multiple testing through such Bayesian hierarchical modeling.  See \cite{bogdan2011asymptotic} and \cite{bogdan2008comparison} for examples of Bayesian multiple testing rules derived as optimal rules with respect to an additive loss functions which are further deiscussed in Section 2. For other Bayesian decision theoretic approaches, see, e.g. \cite{muller2004optimal} , \cite{sun2009large} and \cite{Ghosh2017}\par

The above are examples of multiple hypothesis testing procedures under the assumption of independence of the test statistics for the individual tests. But in practice test statistics may often be dependent. 
It has been observed that when the procedures intended for the independent setup are applied unaltered under dependence, lot of undesirable things can creep in and the performances of these procedures greatly suffer. See in this context, e.g., \cite{qiu2006assessing} and references of Cohen and Sackrowitz for further details. Although these issues have been raised, they have not been adequately resolved in the literature and the area of multiple testing under dependence is still very open to say the least.
\par

The above works on dependence,  
do not focus on the decision theoretic aspect of multiple testing. This aspect has been largely ignored except some references like \cite{sun2009large}, \cite{xie2011optimal}, \cite{cohen2005characterization}, \cite{cohen2007more} and \cite{cohen2008multiple}. In \cite{sun2009large} a decision theoretic study was carried out when the unknown parameters are assumed to be random in nature with a Markovian dependence structure among them. Under the dependent set up, in multiple hypothesis testing procedure there are some methods for estimating FDP. See e.g \cite{fan2017estimation}, \cite{fan2017farm} etc. In these methods it is assumed that the dependency comes into play in the form of some common factors. These methods perform well in presence of factor type dependence set up. There are some other methods by Efron, (e.g. \cite{efron2007correlation}, \cite{efron2010correlated}), where the z-scores are transformed into count data. This translates the problem to the estimation of distribution of the correlations. But still these methods are very problem specific and mostly perform under the assumption of sparsity.
\par 

A natural question is what would the optimal rule (Bayes Oracle) look like under arbitrary form of dependence among test statistics and what would be its asymptotic risk properties in the asymptotic framework of \cite{bogdan2011asymptotic} and \cite{bogdan2008comparison}. If one comes to think of it carefully, this in itself is a very challenging problem under dependence and the reason will become clear shortly. This was left as an open problem in \citep{bogdan2011asymptotic}. Ours is a modest attempt, to work in the direction of this challenging problem. We restrict ourselves to the set-up where the test statistics jointly have a mixture multivariate normal distribution in permutationally invariant setup. We further assume that the parameters have jointly a multivariate normal distribution, given values of a vector of independent Bernoulli random variables. (Details are given in Section 2). With respect to additive loss, the general form of the Bayes Oracle is very easy to derive even under dependence but it is intractable to work with. \par 
We tried to view the problem in  different light and proposed some methods for solving the problem under the set up described above. An algorithm is proposed which converges to a rule that behaves very close to the ideal classifier. Since the signals and noises are different by their variances, it is intuitive to propose some random classifiers which can produce a threshold for separating observations with difference in variances. The additive loss function is tractable in this case through some standard approximations and which in turn makes the problem of risk minimization more manageable. Restricting to constant threshold the risk becomes a continuous function of the classifier through this approximation and hence it can be minimized easily by simple mathematics.\par 
 It is very hard to actually evaluate the performance of these methods because there is no standard popular method in the literature for this kind of a situation. In order to overcome that we need a classifier that can be used as a reference i.e. lower bound on expected error, to be precise. But since this is a dependent set up and it is hard to find a simple "Oracle" type rule due to intractability of the naive rule in \ref{SimRule}, we need to devise an ideal oracle classifier. This classifier, which minimizes the total error i.e. the sum of false positives and false negatives, is then calculated by grid search technique using brute force. This experiment is repeated a significant number of times and averaged over to get the expected 'Optimal' total error of the whole process. It is optimal in a sense that it provides a lower bound on the error rate on all possible classifiers uniformly. Though some of our methods perform close to the ideal classifier, we should note that ideal risk, i.e. risk corresponding to the ideal classifier (obtained by grid search) is a non-achievable lower error bound even in the limit.

\section{Description of the Problem}

\cite{bogdan2011asymptotic} has a set up with independent normal observations and independent normal prior. In this context they tried to solve the testing problem, with which we start with. But we have gone beyond the set up as to extend the problem in the case of correlated normal set up. Suppose we have $m$ observations $X_1,X_2,...,X_m$ such that the vector $\bm{X} =(X_1,X_2,...,X_m)^{\prime}$ follows a $N(\bm{\mu},\sigma_\varepsilon^2 \Sigma)$ distribution given $\bm{\mu} = (\mu_1,\ldots,\mu_m)^{\prime}$. Here $\bm{\mu}$ (unknown) represents the $(m \times 1)$ vector of effects under investigation and $\sigma_\varepsilon^2 \Sigma_{m \times m}$ represents the variability of the random noise (e.g. measurement error). Here we assume that $\sigma_\varepsilon^2 > 0$ is known and $\Sigma$ is a symmetric matrix with $\rho$'s in off diagonal and $1$'s in diagonal, where $\rho$ is known. Thus, given $\bm{\mu}$, $\bm{X}$ is multivariate normal with the correlation between each pair of its components being $\rho$. We further assume that $\bm{\mu}$ is a random vector whose distribution is determined by the values of $m$ unobservable independent Bernoulli($p$) random variables $\nu_i$, for some $p\in (0,1)$. We call $H_{0i}$ the event that $\nu_i=0$ while $H_{1i}$ denotes the event $\nu_i =1$. When $H_{0i}$ is true, $\mu_i$ is assumed to have a $N(0,\sigma_0^2)$ distribution while under $H_{1i}$, it is assumed to follow a $N(0,\sigma_0^2+\tau^2)$ distribution (where $\tau>0$). We further assume that, given $\bm{\nu}=\bm{\nu_0} = (\nu_{01},\nu_{02},\dots \nu_{0m})^{'}$, the different components of $\bm{\mu}$ are correlated, with the common correlation being the same as that between the components of $\bm{X}$, i.e., $\rho$. The assumption of the same correlation between components of the mean vector and the observation vector seems quite natural. Summing up, given $\bm{\nu}=\bm{\nu_0}$, distribution of $\bm{\mu} $ is the following:
  $$
  \bm{\mu}|(\bm{\nu}=\bm{\nu_0})\sim N(0,D_{\nu_0}\Sigma D_{\nu_0}),
  $$
  where $D_{\nu_0}$ is a diagonal matrix with $(D_{\nu_0})_{ii}=\sigma_0 $ if $\nu_{0i}=0$ and $ \sqrt{(\sigma_0 ^2+\tau^2)}$ if $\nu_{0i}=1$ and $\Sigma$ is as described above. When $\sigma_0^2$ is zero or very close to zero and $\tau^2$ is substantially large compared to $\sigma_{0}^2$, then $H_{0i}$'s correspond to the small/insignificant signals or noises while the $H_{1i}$'s correspond to the large or important signals. Our problem is to simultaneously test which components of $\bm{\mu}$ correspond to the large signals, i.e to test $H_{0i}: \nu_i=0$ versus $H_{1i}: \nu_i=1$ simultaneously for $i=1,\ldots,m$. It may be observed that when $\sigma_0^2=0$, $H_{0i}$ corresponds to the point null hypothesis that $\mu_i=0$. Since under $H_{1i}$, $P(\mu_i=0)=0$, our testing problem is equivalent (when $\sigma_0=0$) to the canonical testing problem of $\mu_i=0$ versus $\mu_i \neq 0$. Now define $p_{\bm{\nu_0}}$ as the probability of $\bm{\nu}$ taking the value $\bm{\nu_0}$ i.e. $p_{\bm{\nu_0}}= p^{||\bm{\nu_0}||}(1-p)^{m-||\bm{\nu_0}||}$ with $||\bm{\nu_0}||=$number of 1's in $\bm{\nu_0}$ vector.
 Thus the marginal distribution of $\bm{X}$ is as follows:
  \begin{align} \label{marg1}
  \bm{X}\sim \sum_{\bm{\nu_0}} p_{\bm{\nu_0}} N(0, \sigma_{\varepsilon}^2\Sigma+D_{\bm{\nu_0}}\Sigma D_{\bm{\nu_0}}).
    \end{align} From the above calculation we can see that the marginal distribution is not equi-correlated in the set up used by us. But at the same time we can see that the permutation invariance property is still valid in this case i.e the decision problem is still permutation invariant. $v_i=1$ is concluded, i.e. null hypothesis is rejected if we can see that $X_i$ is greater than some symmetric function of the data set (this includes the case when the function is constant).

In our quest for finding a good multiple testing procedure in our problem, we would ideally look for a procedure that does a decent job of correctly identifying the signals (big signals) from the noises (insignificant signals) by reducing the expected loss.
 Our chosen loss function is an additive one that defines the overall loss in a multiple testing procedure as the sum of losses incurred in the individual testing problems. 
The simplest loss of this kind is the sum of the total number of type I and type II errors made by a multiple testing rule. This loss was originally proposed in Lehmann (\cite{lehmann1957theory}, \cite{lehmann1957theory2}) and later considered by many others. See in this context \cite{sun2007oracle}, \cite{bogdan2011asymptotic}, \cite{bogdan2008comparison}, \cite{datta2013asymptotic} and \cite{sun2009large}. We say that a loss of $\delta_{0}$ is incurred for the $i$-th testing problem when $H_{0i}$ is true but it is rejected, i.e, an error of type I is made. Whereas, a loss of $\delta_{1}$ is said to be incurred for the $i$-th problem when an error of type II occurs in that problem. Here $\delta_0$ and $\delta_1$ are strictly positive real numbers and possibly dependent on $m$. The overall loss of a multiple testing procedure can then be expressed as $L(\bm{\nu}(\bm{X}), \bm{\nu})=\sum_{i=1}^m \delta_i(\nu_i(\bm{X})-\nu_{i})^2$, where $\bm{\nu}$ denotes the true value of nature and $\bm{\nu}(\bm{X})= (\nu_1(\bm{X}),\ldots,\nu_m(\bm{X}))^{\prime}$ represents the corresponding (random) binary vector of 0's and 1's indicating the decisions obtained from a multiple testing procedure. More precisely, $\nu_i(\bm{X})=0$ if the multiple testing rule accepts $H_{0i}$ and $\nu_i(\bm{X})=1$ if $H_{0i}$ is rejected. Thus in the above, $\delta_i= \delta_{0}$ when $\nu_i=0$ but $\nu_i(\bm{X})=1$, whereas $\delta_i= \delta_{1}$ when $\nu_i=1$ but $\nu_i(\bm{X})=0$. The Bayes risk is defined as $R_{m}=E[L(\bm{\nu}(\bm{X}),\bm{\nu})]$, where $E$ denotes expectation with respect to the joint distribution of $(\bm{X}, \bm{\nu})$. It follows easily that
  \begin{align}  
  R_m &=E^{\bm{\nu}} E[L(\bm{\nu}(\bm{X}),\bm{\nu})|\bm{\nu}]\nonumber\\   \label{Risk1} 
  &=\sum_{i=1}^m\left[\delta_{0}(1-p)t_{1i}+\delta_{1} p t_{2i}\right]. 
  \end{align}
In the above, $t_{1i}$ and $t_{2i}$ denote the probabilities of type I and type II errors incurred for the $i$-th testing problem. \par

 It may be noted that in our setup, the parameter space, the marginal distribution and conditional distribution of $(X_1,X_2,\dots X_m)^T$ remain invariant with respect to permutations, which tells us to consider permutation invariant tests. This immediately implies that $t_{1i}=t_{1}$ and $t_{2i}=t_{2}$ $\forall i$  Applying this, the risk becomes $$ R(\bm{\nu},\bm{\nu^*})=m\left[\delta_0(1-p)t_{1}+\delta_1 p t_{2}\right]$$ Our goal would be to minimize $[\delta_0(1-p)t_{1}+\delta_1 p t_{2}]$ among permutations invariant tests to obtain a good approximate rule in this case.\par

It is easy to see that for this additive loss function, the optimal multiple testing rule is the one which simply applies the Bayes rule (with respect to the given $\delta_0$, $\delta_1$ losses) for each individual test and is given by
\begin{align}
\label{SimRule}
\mbox{Reject }H_{0i} \mbox{ if }\frac{f(\bm{X}|\nu_i=1)}{f(\bm{X}|\nu_i=0)} &> \frac{(1-p)}{p} \frac{\delta_0}{\delta_1},\mbox{ and accept it otherwise}
\end{align}
for each $i=1,\dots,n$, where $f(\bm{X}|\nu_i=j)$ is the marginal density of $\bm{X}$ where $\nu_i=j$ for $j=0,1$. But this rule is mathematically extremely intractable and not implementable in practice and it is almost impossible to find an expression for the type I and type II error rates even asymptotically. This is so since $f(\bm{X}|\nu_i=0)$ and $f(\bm{X}|\nu_i=1)$ are given from the mixture densities $\sum \limits_{\bm{\nu}: \nu_i=0} f(\bm{X}|\bm{\nu}) p_{\bm{\nu}}$ and $\sum \limits_{\bm{\nu}: \nu_i=1} f(\bm{X}|\bm{\nu}) p_{\bm{\nu}}$ respectively, where $f(\bm{X}|\bm{\nu})$ is the density of a Multivariate Normal distribution with mean vector zero and covariance matrix $\sigma^2_{\epsilon}\Sigma+D_{\nu}\Sigma D_{\nu}$.  The naive approach of trying to get directly the optimal rule will get us nowhere.

Note that basically our problem is to identify the coordinates $i \in \{1,\ldots,m\}$ for which $\nu_i=1$ and for which $\nu_i=0$. For both $\nu_i=0$ and $\nu_i=1$, the $X_i$'s come from normal distributions with zero means, but different variances, namely $\sigma_{\epsilon}^2+\sigma_0^2$ and $\sigma_{\epsilon}^2+\sigma_0^2+\tau^2$. The important thing to observe is that the variance under $\nu_i=1$ is much larger provided $\tau^2$ is much larger compared to $\sigma^2_0$ and $\sigma_{\epsilon}^2$. So an intuitively appealing approach for finding coordinates with $\nu_i=1$ would be look for integers between 1 to $m$ for which the observed data points appear to come from a distribution with higher variance. If we can find observations which indicate higher variance, the corresponding null hypotheses are more likely to be false. The following lemma, to be proved in the appendix, would indicate that this amounts to looking for $X_i$'s such that $X_i^2$'s are the largest.

 \begin{lemma} \label{lem12}
  \begin{itemize}\label{lemma1}
  \item[(a)] Let  $(X_1,X_2,\dots X_m)^{'} $ follow a Multivariate Normal distribution with correlation matrix $\mathcal{R}$ and their variances are $\sigma_1^2,\sigma_2^2,\dots\sigma_m^2$ respectively. Then $X_i^2\leq_{st} X_j^2$ if and only if $\sigma_i^2\leq \sigma_j^2$.
  \item[(b)]\label{lemma2} Under the assumption of part (a), with equi-correlated correlation matrix $\mathcal{R}$, $X_i^2|Z\leq_{st}X_j^2|Z$ if and only if $\sigma_i^2\leq\sigma_j^2$ where $Z$ is a subset of $\{X_1,X_2,\dots X_m\}$ not containing $X_i$ and $X_j$. 
  \item[(c)] \label{lemma3} Let $X\sim N(0,\sigma_{\epsilon}^2\Sigma+D\Sigma D)$ where $\Sigma$ is equi-correlated and $D=\text{diag}\left(\sigma_1^2,\sigma_2^2,\dots ,\sigma_m^2\right)$. Then $X^2_i|Z\leq_{st} X_j^2|Z$ if and only if $\sigma_i^2\leq \sigma_j^2$ where $Z$ is a subset of $\{X_1,X_2,\dots X_m\}$ not containing $X_i$ and $X_j$. 
\end{itemize}
\end{lemma}
\textbf{Note:} In the above the notation $X \leq_{st} Y$ means $Y$ is stochastically larger than $X$.\par 
Now obviously one has to choose highest $K$ order statistics of the $X_i^2$ in order to execute the above idea. The question is how to choose this $K$ in an 'optimum' way. Put another way, this is equivalent to finding an 'optimum' threshold $C$ (fixed or data dependent) such that we reject $H_{0i}$ whenever $|X_i| > C$. This will precisely be discussed in the next section. Under the restriction to permutation invariant rules, if $C$ is random, it must be a symmetric function of $X_1,X_2,\dots,X_m$.
   
\section{Various Approaches}
In the previous section it is described that the problem of identifying the signals from noises in this context is  dependent on the measure of variance of the coordinates of the random variable. So in a way, it is a clustering problem where we know that there are two clusters. This reduces the problem in finding an appropriate classifier. In general the risk minimization is considered as an optimum way for deciding a rule for testing. \par
 The main idea for the solution of this problem is to classify $X_1,X_2,\dots X_{m}$ into two groups with a classifier of the form $X^2_i>C$ for some $C$ ($C$ may be random). We have given several methods for classifying the variables in two groups. The methods will be described in the next few paragraphs.
 
 \subsection{Algorithm}
 From the previous sections, we can see that the problem is one dimensional and the test statistics have normal distributions with large or small variances. Our goal is to find out an algorithm that gradually converges to some targeted rule. Let us define three quantities first.
 $$T_1 = \left(\frac{\sum_{i=1}^m | X_i |}{m}\right)$$
 $$T_2 = \left(\frac{\sum_{i=1}^m X_i^2}{m}\right)^\frac{1}{2}$$
 $$T_4 = \left(\frac{\sum_{i=1}^m X_i^4}{m}\right)^\frac{1}{4}$$
We have proposed an iterative algorithm for determining the classifier C, which works good with reasonable error of false positive and false negative. To find the C we do the iterative steps in the following way:
       \begin{itemize}
       \item[1.] {\bf Initialization:} Start with $Z_0=T_i$ for some $i\in \{1,2,4\}$. Classify the vector of co-ordinate wise absolute value of $X$ with this classifier.
       \item[2.] {\bf Loop:} The co-ordinates of $\bm{X}$ for which the corresponding absolute values are less than $Z_0$ and those which are greater than $Z_0$ form two groups of absolute values. Call the group means $A_1$ and $A_2$ respectively and obtain $Z_1=\frac{A_1+A_2}{2}$.
       \item[3.] Go to step 1 with $Z_1$ and obtain $Z_2,Z_3,\dots$ respectively. 
       \item[4.] {\bf Termination:} Terminate the process in the i-th step if $|Z_{i+1}-Z_i|<f$, where $f$ is a predetermined very small tolerance value. 
       \end{itemize}
       This is because it can be shown that the sequence $\{Z_n\}_{n\geq 0}$ converges in general.
From the algorithm, we get a reasonable $C$ with iterative limit of $Z_i$'s. Different points for initializing the algorithm were studied and the performance are reported in the simulated data set. While selecting different statistics as points of initialization, homogeneity among the co-ordinates are maintained because observations only differ in variability. We can also initialize the algorithm at the grand mean or a grand median point of the data set. The following result shows the rationale behind proposing this algorithm.
 \begin{result} \label{Algo}
 Let $w_1,w_2,\dots,w_m$ be $m$ positive observations. Then within group variance:$$V_w(C)=\sum_{w_i\leq C}\left(w_i-\bar{w_1}\right)^2+\sum_{w_i> C}\left(w_i-\bar{w_2}\right)^2$$ is minimized for a value of $C$ which necessarily satisfies $C=\frac{\bar{w_1}+\bar{w_2}}{2}$ with $\bar{w_i}$ giving the mean of the $i$-th group.
 \end{result}  

The algorithm is an iterative one i.e. in each step of the algorithm, the within group variance is decreasing and later the data is classified based on that classifier. Therefore in each step the classifier increasingly reduces the within group variance and hence forces the classifier to divide the data into two clusters. This phenomena also leads to less expected mis-classification.  The algorithm is similar to many well-known algorithms which converges for a much wide region of non-pathological data sets.
 
\subsection{Proposed Random Classifier}
 In the previous section, we have devised an algorithm to lead us to a random classifier (as $Z_i$'s are functions of the data sets), which will do a decent job in approximating the optimal test. Various points of initialization (though similar in output) were used to see the performance and convergence of the process. We have already stated that the rule for classifying the observations in two groups will involve a classifier, which can be a constant as well as random. The statistics defined in the previous section i.e. $T_1, T_2, T_4$ might be a good random choice of classifier for the data points that differ in variability. These classifiers fall in the category of random classifier and their performances are checked in simulation study in the next section.\par 
 
 Let us define the following general quantity : $$T_{2h}=\left(\frac{\sum_{i=1}^m X^{2h}_i}{m}\right)^\frac{1}{2h}$$ Since the main interest of the observations are their measure of variability, any even powered moment is a potential choice for the random classifier. The $2h$-th root is taken to maintain homogeneity among the sequence of statistics.\par 
 
  Dependent on the data and its labels we have calculated the optimal values of $h$, for which corresponding $T_{2h}$ performs well. We have given a table for the corresponding optimal rules of $h$ in Appendix B.\par 
  In what follows, we have found optimal value of $C$ which works as a constant classifier to minimize risk, which do not depend on the data and its labels. In a similar way we may find an optimal value of $h$, which works for a good classifier depending only on $(\sigma^2_0,\sigma^2_{\epsilon},\tau^2,m,\rho,\beta)$. We expect that will perform close to the optimal classifier or ideal classifier as can be seen from  table \ref{Err-CF-Opt}. It is a worthwhile task to find $h_{opt}$ as discussed above.

\subsection{Risk Minimization C}
In this section we will first provide a rule,close to the optimal rule under the additive loss. As mentioned in Section 2, our search for the optimal rule will be restricted within the class of permutation invariant rules. To start with, consider a multiple testing rule that rejects $H_{0i}$ whenever $\{|X_i| > C \}$ for $i=1,\ldots,m$, where $C$ is a fixed constant. The expected loss of this rule is given by

    \begin{align}
    R(C)=\delta_0m(1-p)P\left[|Y_1|>\frac{C}{\sqrt{\sigma_{\epsilon}^2+\sigma_0^2}}\right]+\delta_Amp P\left[|Y_1|<\frac{C}{\sqrt{\sigma_{\epsilon}^2+\sigma_0^2+\tau^2}}\right],
     \end{align}
where $Y_1$ denotes a standard normal random variable. This follows easily from the model assumed and observing that for the multiple testing rule under consideration, the errors are identical for each $i$. When $X_i$'s and $\mu_i$'s are independent for $i=1,2,\dots,m$, the optimal test has a rejection region of this kind (with the threshold independent of the observations but dependent on the model parameters). This rule is called a Bayes Oracle since it depends on the unknown model parameters. This has been explicitly derived in \cite{bogdan2011asymptotic} (see equation 2.4 of Section 2 of the paper). They also calculated, under their asymptotic framework, the asymptotic risk of the Bayes Oracle. \par

As mentioned in the previous section, we would ideally like to fully extend the study of optimality in our dependent setup. But as already explained, the optimal test (the Bayes Oracle) in our dependent setup does not come out as a simple thresholding rule as in the independent case. Instead, it is analytically very complicated, making the corresponding risk quite intractable. However, it is curious to note that the risk function for the fixed threshold test above is the same for both the dependent and the independent cases if the marginal distributions are the same, since the risk depends only on the marginal distributions of the observations. So one can easily minimize the risk with respect to $C$ to get a critical region that minimizes the risk among fixed threshold rules even for the dependent case. Note that the optimal fixed threshold rule for independent case is also optimal among fixed threshold rules in the dependent case. In dependent cases the general optimal rule beats the above rule.

As mentioned before, here we set ourselves a moderate goal to get to this optimal procedure with fixed threshold. For simplicity we will only consider the case when $\delta_0=\delta_A = 1$. \par 

Our goal now is to try to derive its asymptotic risk. We restrict ourselves to providing a heuristic approximation to the exact asymptotic or fixed sample risk. Towards that, we recall the risk of a fixed threshold rule (with fixed threshold $C$) described at the beginning of the section. The expression of $R$ can be approximated in the following way. Assuming $\frac{C}{\sqrt{\sigma_{\varepsilon}^2+\sigma_0^2}}$ is large and $\frac{C}{\sqrt{\sigma_{\varepsilon}^2+\sigma_0^2+ \tau^2}}$ is small, we have
     \begin{align}
     R(C)& \approx \delta_0 m (1-p) \frac{\sqrt{\sigma_0^2+\sigma_{\varepsilon}^2}\sqrt{2}}{C\sqrt{\pi}} e^{-\frac{C^2}{2{(\sigma_0^2+\sigma_{\varepsilon}^2)}}} + \delta_A m p \frac{C\sqrt{2}}{\sqrt{\pi(\sigma_0^2+\sigma_{\varepsilon}^2+\tau^2)}} \nonumber\\
     &= \frac{V}{C}e^{-\frac{C^2}{2{(\sigma_0^2+\sigma_{\varepsilon}^2)}}}+UC \nonumber\\
     &=f(C)
     \end{align}
     where $V= \delta_0 m (1-p)\sqrt{\sigma_0^2+\sigma_{\varepsilon}^2}\sqrt{\frac{2}{\pi}}$ and $U=\frac{\delta_A m p}{\sqrt{(\sigma_0^2+\sigma_{\varepsilon}^2+\tau^2)}} \sqrt{\frac{2}{\pi}} $. For the first summand, we exploit the fact that $\frac{C}{\sqrt{\sigma_{\varepsilon}^2+\sigma_0^2}}$ is large and employ the standard approximation to normal tails using Mill's Ratio. For the second summand, we note that since $\frac{C}{\sqrt{(\sigma_{\varepsilon}^2+ \sigma_0^2+ \tau^2)}}$ is small and for small $x$, $P[|N(0,1)|<x]\approx 2x\phi(0)$.
     Now observe that $f(C)$ as defined above is a convex function of $C$, since 
     \begin{align*}
     f'(C)&=U-\frac{V}{C^2}e^{-aC^2}-2aVe^{-aC^2}\\
          &=U-Ve^{-aC^2}\left(\frac{1}{C^2}+2a\right)  
      \end{align*}
where $a=\frac{1}{2(\sigma_0^2+\sigma_{\varepsilon}^2)}$ and $f'(C)$ is an increasing function of $C$. {as $f^{\prime \prime}(C) > 0$.} 
       
\begin{align*}
f^{'}(C) &= U - 2aVe^{-aC^2} -\mathcal{O}\left(\frac{1}{e^{aC^2}C^2}\right)
\end{align*} 
By ignoring the third term in the expression of $f^{'}(C)$, we get the maximizer $C$, which is approximately given by the following:
\begin{align*}
       U&=2aVe^{-aC^2}\\
       C&=\sqrt{\frac{1}{a}log\left(\frac{2aV}{U}\right)}
       \end{align*}
\begin{align}
\text{Putting the values:}\nonumber\\ \label{C} 
C&\approx \sqrt{2(\sigma_0^2+\sigma_{\varepsilon}^2)\log\left(\frac{\delta_0 (1-p)}{\delta_A p}\sqrt{\left(1+\frac{\tau^2}{\sigma_0^2+\sigma_{\varepsilon}^2}\right)}\right)}
 \end{align}

\subsection{Ideal Classifier}
 In the previous subsections various methods for attacking the problem of multiple hypothesis testing under the difference in variability set up is described. In this context the natural thing to ask is : How these methods are performing? Since we do not have any standard oracle rule, it is hard to evaluate the credibility of the methods. An ideal rule, which can be used as a reference frame while calculating the degree of accuracy of these methods in action, is required.\par 
 
 While devising the methods for our problem, we needed an ideal choice of $C$, so that we could compare our methods with the ideal case and get an estimate of efficiency. Now the problem lies in selection of the ideal $C$ i.e. the best performing $C$ while classifying the simulated data. In order to choose that particular $C$ grid search technique is applied to approximately find the closest ideal $C$ to work with. It is observed that we have to classify according to a measure of variance. Therefore, we have taken modulus value of the simulated $X$ vector coordinate-wise. Now we looked for the range of the absolute values of the co-ordinates of $X$ and started with the lowest point of the range as our starting $C$.  After classification with this $C$ and finding out the error (Sum of false positive and false negative), an increment value for $C$ is fixed. Every time $C$ is increased with the corresponding fixed increment value, the simulated data (from a pre-decided fixed set up, i.e. fixing $m,\sigma_{\epsilon},\sigma_0,\tau,$ and $\rho$) is classified and the total error in the case is calculated. In this way we can get a sequence of total errors for different values of $C$. Among them, the $C$ which corresponds to the minimum total error among these sequence of errors is selected as our ideal $C$ and the corresponding total error is reported as ideal case total error.\par

 The process of selecting the best $C$ and determining the ideal case total error by brute force is discussed in the last paragraph. But there are some subtle points which needs to be understood properly in order to justify the process. As the process goes on finding the classifier, there may be multiple choices of best $C$, which provides the same total error at the end. But we do not need to bother about those multiple $C$'s in this case at all, because our purpose for generating this ideal case is to get the idea of the total error in the best possible choice of $C$ and create an ideal case to compare with for testing our proposed methods for choice of $C$. So any one of these best $C$'s can be selected as our ideal $C$ and the simulated data will be classified to get the total error.\par

\begin{rem}
Note that, we have made classification in the ideal case with the knowledge of which observation comes from which $\sigma_i^2$. Here in practice, it may not be achievable as the classifier is not a function of $Y_1^2,Y_2^2,\dots Y_{m}^2$ alone. Thus the ideal case can be looked upon as some lower bound which may not be achievable even in the limit.
\end{rem}

\section{Simulation and Discussion}
As we have stated earlier, our problem now boils down to finding a suitable $C$ which will classify the observations coming from an underlying set up. We did not have any real standard data to test run our process, so we have simulated and performed the tests to check the validity of our method and compared the cases having non-zero correlation terms with the independent cases. As we have developed our method, we have shown that our process with non-zero correlation coefficient in the equi-correlated set up is at least as good as the independent case in terms of the risk function.\par

In the following paragraph, the observations obtained directly from the simulation studies where the correlation coefficients are equal and non-negative will be described. It is stated earlier that $C$ is a symmetric function of $X_1, X_2, \dots X_m$. For this reason, some typical functions are selected to study the behavior of the total error. From the simulation studies, we have some good observations to make about our procedure. We claimed that our procedure is at least as good as in the case of independent set up in terms of risk function. This phenomena is reflected in the simulation studies as well. If we go through the table we can easily see that irrespective of the method of choosing the $C$, our claim holds. Apart from that as the correlation gets high, these methods perform even better. In the process  we have assumed that the $\nu_i$'s follow $Ber(p)$. We have used the same setting for choosing $p$ as in \cite{bogdan2011asymptotic}. $p$ is selected as $m^{-\beta}$ while running the simulation. The simulation for each set of parameter value is run 10000 times and then the average of the sum of false positive and false negative is taken to get an estimate of the expected misclassification in each case. \par 

There is another observation that can be made from the simulated data sets. In this process $\tau$ is generally assumed to be larger than $\sigma_0$, which in turn helps the process of classification in this way to work properly. This is because we are dividing the co-ordinates of the simulated random variable into two parts according to a measure of variance and assigning the $\nu_i$ to be 1 for having larger variance and 0 to the rest of them. Now as we can see from the simulation studies that if the ratio $\frac{\tau}{\sigma_0}$ is large then the classification is good and the expected number of both false positive and false negative decrease which is expected. $T_4$ performs better compared to $T_1$ and $T_2$ if the ratio is high.\par

In case of C determined by the expression (\ref{C}) stated above, we can see that the expression depends on $\delta_0$ and $\delta_A$. In our set up we have assumed that these two values are same and equal to 1 and carried out the simulation with C as in expression (\ref{C}). From the simulation table of total error presented here, we can see that this choice of $C$ works better if $\beta$ is higher i.e. it is sensitive to $\beta$ and it gets extremely close to the ideal case
for higher values of $\beta$. In all the cases, performance of $T_1$ and $T_2$ are more or less similar to each other. The algorithm performs better than the two expressions $T_1$ and $T_2$, when the ratio $\frac{\tau}{\sigma_0}$ is not comparatively low. But the algorithm works better than the two expressions mentioned above as $m$ and $\beta$ becomes larger over all. Another interesting observation is that for $T_1$,$T_2$ and for the choice of $C$ with the algorithm stated above, we can see that for smaller value of $\beta$ the bias is toward the false negative values i.e. expected false negative is almost uniformly larger than the expected false positives. But in case of higher values of $\beta$ the bias gets reversed. In case of
$C$ determined by the formula derived above(in the last section) the bias is almost always toward the false negative cases. As $m$ increases all the possible choices of $C$ works better gradually. \par

 Now to compare the performance of the methods we may calculate the discrepancies of various methods using the following formula:
       $$\textbf{Discripancy in Percentage}=100\times\frac{E_K-E_{K_0}}{E_K}$$
       where $E_K$ is the total error in the corresponding choice of $C$ and $E_{K_0}$ is the total error in the ideal choice of $C$. Here by total error we mean the sum of the expected number of false positive and expected number of false negative cases.

\section*{Acknowledgment }
   We are grateful to Professor Malay Ghosh for his valuable suggestions and comments. \par 

\bibliographystyle{plainnat}
\bibliography{Ref}

\newpage
\section*{Appendix A}
\begin{enumerate}
\item Proof of Equation (\ref{Risk1}):
   \begin{align*}   
  R(\bm{\nu},\bm{\nu^*})&=EE[L(\bm{\nu},\bm{\nu^*})|\bm{\nu}=\bm{\nu_0}]\nonumber\\ 
  &=\sum_{\nu_0} p_{\nu_0} E\left[\sum_{i=1}^m \delta_0 {1}_{[(\nu_{0i}-\nu_i^*)=1]}+\sum_{i=1}^m \delta_A {1}_{[(\nu_{0i}-\nu_i^*)=-1]}\right]\nonumber\\
  &=\sum_{i=1}^m E\left[\sum_{\bm{\nu_0}} p_{\nu_0}\left(\delta_0 {1}_{[(\nu_{0i}-\nu_i^*)=1}+\delta_A {1}_{[(\nu_{0i}-\nu_i^*)=-1]}\right)\right]\\
  &=\sum_{i=1}^m \left[\sum_{\bm{\nu_0}} p_{\nu_0}\left(\delta_0 E({1}_{(\nu_{0i}=1)}|\nu_i^*=0)P(\nu_i^*=0)+\delta_A E({1}_{(\nu_{0i}=0)}|\nu_i^*=1)P(\nu_i^*=1)\right)\right]\\
  &=\sum_{i=1}^m\delta_0(1-p)t_{1i}+\delta_A p t_{2i} 
  \end{align*}

\item Proof of Lemma(\ref{lem12}):
   \begin{itemize} 
   \item[(a)] First we consider marginals of $X_i^2$ and $X_j^2$ with their marginal variance $\sigma_i^2$ and $\sigma_j^2$ as their variances respectively. Then 
 $$X_i^2\sim \sigma_i^2 V$$, where $V$ is chi square with one degree of freedom. From this it easily follows that $X_i^2\leq_{st}X_j^2$ if and only if $\sigma_i^2\leq\sigma_j^2$.    
   \item[(b)] Here we shall prove $X_i^2|Z\leq_{st}X_j^2|Z$ when the inequality in variance holds as stated above. Here $Z$ is a subset of $\left\{X_1,X_2,\dots X_m\right\}$ deleted by $X_i$ and $X_j$ respectively. Now in the equi-correlated set up, without loss of generality, instead of $i\neq j$ we may simply work with $1$ and $2$. Define
  $$U=\left[\left(\frac{X_1}{\sigma_1},\frac{X_2}{\sigma_2}\right)\,\middle|\,\left(\frac{X_3}{\sigma_3},\frac{X_4}{\sigma_4}\dots \frac{X_m}{\sigma_m}\right)\right]$$ 
  This quantity is free of $\sigma_i^2$ and $\left(\frac{X_1}{\sigma_1},\frac{X_2}{\sigma_2},\frac{X_3}{\sigma_3},\frac{X_4}{\sigma_4}\dots \frac{X_m}{\sigma_m}\right)$ has exchangeable distributions. Now $U=(U_1,U_2)$ which are exchangeable. $(X^2_1,X^2_2)=(\sigma^2_1U^2_1,\sigma^2_2U^2_2)$ which has equi-correlated matrix
  $\mathcal{R^*}$. Hence by part (a) the result follows.
  \item[(c)] This part follows from part (a) and part (b).
  \end{itemize}

\item Proof of Result {\ref{Algo}}:\\
   Let us assume a continuous p.d.f. of $w$ and call it $f_w$. Then within group variance $V_W(C)$ of the two groups obtained from $w$, using $C$, is a continuous function of $C$. We consider $$\frac{\partial V_W(C)}{\partial C}=0$$ and obtain the result. Now as the result holds for continuous p.d.f., it is easy to see that it holds for the discrete case also.
\end{enumerate}

\newpage
\section*{Appendix B}

{\setlength\tabcolsep{3.5pt}
\begin{longtable}[c]{|c|c|c|c|c|c|c|c|c|c|c|c|c|}
\caption{Comparison of Total Errors among Clustering Classifiers}
\label{Err-CF-Tab}\\
\hline
$\sigma^2_{\epsilon}$ & $\sigma^2_0$ & $\tau^2$ & $m$  & $\rho$ & $\beta$ & $T^{Algo}_1$ & $T^{Algo}_2$ & $T^{Algo}_4$ & $T_1$  & $T_2$  & $T_4$  & Ideal    \\ \hline
\endfirsthead
\multicolumn{13}{c}%
{{\bfseries Table \thetable\ continued from previous page}} \\
\hline
$\sigma^2_{\epsilon}$ & $\sigma^2_0$ & $\tau^2$ & $m$  & $\rho$ & $\beta$ & $T^{Algo}_1$ & $T^{Algo}_2$ & $T^{Algo}_4$ & $T_1$  & $T_2$  & $T_4$  & Ideal    \\ \hline
\endhead
1                     & 1            & 15       & 100  & 0      & 0.2     & 0.205        & 0.2047       & 0.2119       & 0.2434 & 0.2009 & 0.2242 & 0.1786   \\ \hline
2                     & 3            & 90       & 100  & 0      & 0.2     & 0.1753       & 0.1753       & 0.1882       & 0.1615 & 0.1502 & 0.2109 & 0.1263   \\ \hline
1                     & 1            & 15       & 500  & 0      & 0.2     & 0.1678       & 0.1678       & 0.1687       & 0.2532 & 0.1803 & 0.1775 & 0.1589   \\ \hline
2                     & 3            & 90       & 500  & 0      & 0.2     & 0.1837       & 0.1837       & 0.1867       & 0.1614 & 0.1524 & 0.219  & 0.1378   \\ \hline
1                     & 1            & 15       & 1000 & 0      & 0.2     & 0.1544       & 0.1544       & 0.1544       & 0.26   & 0.1762 & 0.1605 & 0.149    \\ \hline
2                     & 3            & 90       & 1000 & 0      & 0.2     & 0.1381       & 0.1381       & 0.1394       & 0.1743 & 0.1149 & 0.1577 & 0.1111   \\ \hline
1                     & 1            & 15       & 100  & 0.1    & 0.2     & 0.2199       & 0.2198       & 0.2294       & 0.2414 & 0.2113 & 0.2467 & 0.1899   \\ \hline
2                     & 3            & 90       & 100  & 0.1    & 0.2     & 0.2086       & 0.2087       & 0.2238       & 0.1651 & 0.1833 & 0.2578 & 0.1437   \\ \hline
1                     & 1            & 15       & 500  & 0.1    & 0.2     & 0.1756       & 0.1756       & 0.1768       & 0.2495 & 0.1833 & 0.1883 & 0.1658   \\ \hline
2                     & 3            & 90       & 500  & 0.1    & 0.2     & 0.1661       & 0.1661       & 0.169        & 0.1623 & 0.1368 & 0.1958 & 0.1269   \\ \hline
1                     & 1            & 15       & 1000 & 0.1    & 0.2     & 0.1565       & 0.1565       & 0.1566       & 0.2586 & 0.1763 & 0.1633 & 0.1506   \\ \hline
2                     & 3            & 90       & 1000 & 0.1    & 0.2     & 0.1489       & 0.1489       & 0.1503       & 0.1678 & 0.1223 & 0.1723 & 0.1178   \\ \hline
1                     & 1            & 15       & 100  & 0.8    & 0.2     & 0.2138       & 0.2147       & 0.2214       & 0.1941 & 0.21   & 0.2748 & 0.1618   \\ \hline
2                     & 3            & 90       & 100  & 0.8    & 0.2     & 0.1484       & 0.1493       & 0.158        & 0.119  & 0.1372 & 0.208  & 0.0949   \\ \hline
1                     & 1            & 15       & 500  & 0.8    & 0.2     & 0.1467       & 0.146        & 0.1448       & 0.1956 & 0.1463 & 0.1544 & 0.1177   \\ \hline
2                     & 3            & 90       & 500  & 0.8    & 0.2     & 0.1128       & 0.1128       & 0.1142       & 0.1174 & 0.0968 & 0.1424 & 0.08     \\ \hline
1                     & 1            & 15       & 1000 & 0.8    & 0.2     & 0.1401       & 0.1394       & 0.138        & 0.1994 & 0.1427 & 0.1432 & 0.1124   \\ \hline
2                     & 3            & 90       & 1000 & 0.8    & 0.2     & 0.0922       & 0.0921       & 0.0927       & 0.128  & 0.0819 & 0.1101 & 0.0681   \\ \hline
1                     & 1            & 15       & 100  & 0      & 0.9     & 0.3022       & 0.285        & 0.2462       & 0.4041 & 0.2867 & 0.1459 & 0.0129   \\ \hline
2                     & 3            & 90       & 100  & 0      & 0.9     & 0.2166       & 0.2015       & 0.1349       & 0.3746 & 0.2277 & 0.079  & 0.0141   \\ \hline
1                     & 1            & 15       & 500  & 0      & 0.9     & 0.3166       & 0.3086       & 0.3033       & 0.4164 & 0.3036 & 0.1553 & 0.0058   \\ \hline
2                     & 3            & 90       & 500  & 0      & 0.9     & 0.3232       & 0.3167       & 0.3113       & 0.4212 & 0.3092 & 0.1612 & 0.0011   \\ \hline
1                     & 1            & 15       & 1000 & 0      & 0.9     & 0.3294       & 0.322        & 0.3194       & 0.4239 & 0.3156 & 0.1828 & 8.00E-04 \\ \hline
2                     & 3            & 90       & 1000 & 0      & 0.9     & 0.3259       & 0.3204       & 0.3185       & 0.4231 & 0.3131 & 0.1701 & 6.00E-04 \\ \hline
1                     & 1            & 15       & 100  & 0.1    & 0.9     & 0.3272       & 0.3092       & 0.2806       & 0.4167 & 0.3046 & 0.1698 & 0.0065   \\ \hline
2                     & 3            & 90       & 100  & 0.1    & 0.9     & 0.2167       & 0.201        & 0.1342       & 0.375  & 0.2275 & 0.0782 & 0.0138   \\ \hline
1                     & 1            & 15       & 500  & 0.1    & 0.9     & 0.3359       & 0.327        & 0.3233       & 0.4262 & 0.3194 & 0.1916 & 0        \\ \hline
2                     & 3            & 90       & 500  & 0.1    & 0.9     & 0.3181       & 0.3108       & 0.3019       & 0.4191 & 0.3026 & 0.1395 & 0.0021   \\ \hline
1                     & 1            & 15       & 1000 & 0.1    & 0.9     & 0.3347       & 0.3268       & 0.3246       & 0.4261 & 0.3192 & 0.1909 & 0        \\ \hline
2                     & 3            & 90       & 1000 & 0.1    & 0.9     & 0.325        & 0.319        & 0.3167       & 0.4226 & 0.3105 & 0.155  & 0.0011   \\ \hline
1                     & 1            & 15       & 100  & 0.8    & 0.9     & 0.3926       & 0.3652       & 0.3236       & 0.4477 & 0.3436 & 0.1914 & 0.0044   \\ \hline
2                     & 3            & 90       & 100  & 0.8    & 0.9     & 0.1932       & 0.1744       & 0.1159       & 0.3823 & 0.2166 & 0.0665 & 0.0088   \\ \hline
1                     & 1            & 15       & 500  & 0.8    & 0.9     & 0.4099       & 0.3936       & 0.3887       & 0.4579 & 0.3602 & 0.2087 & 0.002    \\ \hline
2                     & 3            & 90       & 500  & 0.8    & 0.9     & 0.3952       & 0.3808       & 0.3411       & 0.4533 & 0.3434 & 0.1483 & 0.0014   \\ \hline
1                     & 1            & 15       & 1000 & 0.8    & 0.9     & 0.4217       & 0.4073       & 0.4061       & 0.4636 & 0.3718 & 0.241  & 5.00E-04 \\ \hline
2                     & 3            & 90       & 1000 & 0.8    & 0.9     & 0.3953       & 0.3829       & 0.3459       & 0.4537 & 0.3436 & 0.1352 & 0.0014   \\ \hline
\end{longtable}
}

\begin{longtable}[c]{|c|c|c|c|c|c|c|c|c|c|}
\caption{Comparison of Total Error among Optimized Classifiers}
\label{Err-CF-Opt}\\
\hline
$\sigma^2_{\epsilon}$ & $\sigma^2_0$ & $\tau^2$ & $m$  & $\rho$ & $\beta$ & $C^{Risk}_{Min}$ & Ideal     & $C^{Moment}_{opt}$ & $h_{opt}$ \\ \hline
\endfirsthead
\multicolumn{10}{c}%
{{\bfseries Table \thetable\ continued from previous page}} \\
\hline
$\sigma^2_{\epsilon}$ & $\sigma^2_0$ & $\tau^2$ & $m$  & $\rho$ & $\beta$ & $C^{Risk}_{Min}$ & Ideal     & $C^{Moment}_{opt}$ & $h_{opt}$ \\ \hline
\endhead
1                     & 1            & 15       & 100  & 0      & 0.2     & 0.204126         & 0.178563  & 0.195695           & 4.4612    \\ \hline
2                     & 3            & 90       & 100  & 0      & 0.2     & 0.14535          & 0.126309  & 0.15009            & 4.0166    \\ \hline
1                     & 1            & 15       & 500  & 0      & 0.2     & 0.1675878        & 0.1589174 & 0.1722134          & 5.088     \\ \hline
2                     & 3            & 90       & 500  & 0      & 0.2     & 0.1441302        & 0.137782  & 0.152417           & 4         \\ \hline
1                     & 1            & 15       & 1000 & 0      & 0.2     & 0.1544564        & 0.1489643 & 0.1599835          & 5.7948    \\ \hline
2                     & 3            & 90       & 1000 & 0      & 0.2     & 0.1149929        & 0.1111424 & 0.1149324          & 4         \\ \hline
1                     & 1            & 15       & 100  & 0.1    & 0.2     & 0.214527         & 0.189941  & 0.208692           & 4.2616    \\ \hline
2                     & 3            & 90       & 100  & 0.1    & 0.2     & 0.163227         & 0.143661  & 0.183251           & 4.0036    \\ \hline
1                     & 1            & 15       & 500  & 0.1    & 0.2     & 0.1748608        & 0.1658236 & 0.1784904          & 4.7182    \\ \hline
2                     & 3            & 90       & 500  & 0.1    & 0.2     & 0.1339354        & 0.1268686 & 0.1368062          & 4         \\ \hline
1                     & 1            & 15       & 1000 & 0.1    & 0.2     & 0.1568779        & 0.1506001 & 0.162047           & 5.6396    \\ \hline
2                     & 3            & 90       & 1000 & 0.1    & 0.2     & 0.1223936        & 0.1177868 & 0.1223492          & 4         \\ \hline
1                     & 1            & 15       & 100  & 0.8    & 0.2     & 0.250247         & 0.161751  & 0.19761            & 3.7784    \\ \hline
2                     & 3            & 90       & 100  & 0.8    & 0.2     & 0.162716         & 0.094927  & 0.136098           & 4.0476    \\ \hline
1                     & 1            & 15       & 500  & 0.8    & 0.2     & 0.1750034        & 0.1177484 & 0.126833           & 4.921     \\ \hline
2                     & 3            & 90       & 500  & 0.8    & 0.2     & 0.1287114        & 0.0799796 & 0.0940992          & 4.1412    \\ \hline
1                     & 1            & 15       & 1000 & 0.8    & 0.2     & 0.1667488        & 0.1123881 & 0.1195987          & 5.0386    \\ \hline
2                     & 3            & 90       & 1000 & 0.8    & 0.2     & 0.1095179        & 0.0680599 & 0.0762898          & 4.2716    \\ \hline
1                     & 1            & 15       & 100  & 0      & 0.9     & 0.015777         & 0.012873  & 0.017634           & 29.3884   \\ \hline
2                     & 3            & 90       & 100  & 0      & 0.9     & 0.01769          & 0.014071  & 0.016089           & 14.8324   \\ \hline
1                     & 1            & 15       & 500  & 0      & 0.9     & 0.006574         & 0.0058144 & 0.0067168          & 35.9596   \\ \hline
2                     & 3            & 90       & 500  & 0      & 0.9     & 0.0013832        & 0.0010694 & 0.0027006          & 45.3904   \\ \hline
1                     & 1            & 15       & 1000 & 0      & 0.9     & 0.0009331        & 0.0007527 & 0.002074           & 63.4846   \\ \hline
2                     & 3            & 90       & 1000 & 0      & 0.9     & 0.0007107        & 0.0005604 & 0.0015375          & 52.4064   \\ \hline
1                     & 1            & 15       & 100  & 0.1    & 0.9     & 0.008606         & 0.006532  & 0.014347           & 37.6474   \\ \hline
2                     & 3            & 90       & 100  & 0.1    & 0.9     & 0.017515         & 0.013781  & 0.015772           & 14.9062   \\ \hline
1                     & 1            & 15       & 500  & 0.1    & 0.9     & 0.0002512        & 0         & 0.0032348          & 67.3786   \\ \hline
2                     & 3            & 90       & 500  & 0.1    & 0.9     & 0.0025996        & 0.0021044 & 0.0029762          & 32.7084   \\ \hline
1                     & 1            & 15       & 1000 & 0.1    & 0.9     & 0.0001331        & 0         & 0.0018177          & 73.8534   \\ \hline
2                     & 3            & 90       & 1000 & 0.1    & 0.9     & 0.0013476        & 0.0011179 & 0.0016696          & 38.439    \\ \hline
1                     & 1            & 15       & 100  & 0.8    & 0.9     & 0.008674         & 0.00438   & 0.010062           & 33.13     \\ \hline
2                     & 3            & 90       & 100  & 0.8    & 0.9     & 0.018041         & 0.008838  & 0.010467           & 12.1456   \\ \hline
1                     & 1            & 15       & 500  & 0.8    & 0.9     & 0.0034084        & 0.0019958 & 0.0031358          & 38.088    \\ \hline
2                     & 3            & 90       & 500  & 0.8    & 0.9     & 0.0026018        & 0.0013576 & 0.0020084          & 24.5464   \\ \hline
1                     & 1            & 15       & 1000 & 0.8    & 0.9     & 0.0008915        & 0.0005178 & 0.001586           & 54.6842   \\ \hline
2                     & 3            & 90       & 1000 & 0.8    & 0.9     & 0.0026875        & 0.0013757 & 0.0015856          & 17.976    \\ \hline
\end{longtable}

\end{document}